\documentclass[12pt]{article}
\usepackage{latexsym}
\usepackage{amssymb}
\usepackage{amsmath}
\usepackage{amsthm}
\usepackage{rotating}
\usepackage{multirow}
\usepackage{mathrsfs}
\usepackage{stmaryrd}
\usepackage{wasysym}
\usepackage{enumerate}
\usepackage{slashed}
\usepackage{nicefrac}
\usepackage{rawfonts}
\input{prepictex}
\input{pictex}
\input{postpictex}
\DeclareMathAlphabet{\zap}{OT1}{pzc}{m}{it}

\newtheorem{main}{Theorem}

\newtheorem{thm}{Theorem}[section]
\newtheorem{prop}[thm]{Proposition}

\newtheorem{cor}[thm]{Corollary}

\def\ZZ{{\mathbb Z}}

\def\CP{\mathbb{CP}}

\begin{document}

\title{Einstein Constants and Smooth Topology, Revisited}

\author{Claude LeBrun and Mingyang Li\\ 
Stony
 Brook
 University} 
  
\date{}

\maketitle

\begin{abstract} 
  In \cite{cat},   Catanese  and  the first author showed   that, for every  $m\geq  2$,   there exists a smooth 
closed  $4m$-manifold that  admits both a $\lambda > 0$ Einstein metric and a 
$\lambda < 0$  Einstein metric; thus,  the sign of the Einstein constant $\lambda$ is 
not always determined by the differential topology of the host compact  manifold. This result was  later improved by  Rasdeaconu and \c{S}uvaina  \cite{raresioana}, who showed that 
there are in fact  at least  $\binom{m + 3}{3}$ topologically distinct such examples  in each dimension $4m$, $m\geq 2$. 
In this note, we further improve  these earlier results by 
showing  that {\em infinitely many} topologically distinct examples of this phenomenon actually  occur  in  
dimension $4m$,  for each $m\geq 2$. 
  \end{abstract}

\section{Introduction}

Recall that a  Riemannian metric $g$ is said to be {\em Einstein} if its Ricci curvature, considered as  a function on the 
unit tangent bundle, is constant. This is of course equivalent to requiring that   the  Ricci tensor $r$ of $g$   satisfy 
\begin{equation}
\label{eineq}
r =\lambda g
\end{equation}
   for some real number $\lambda$, which is then called the {\em Einstein constant} of $g$.  
      An  {\em Einstein manifold}  is then defined to be 
   a compact connected   manifold  (without boundary) that is equipped with a specific Einstein   metric.

  Some of the modern prominence  of the study of \eqref{eineq} 
  can  be attributed to     the  popularity  of 
   the pseudonymous A.\ L.\ 
   Besse's book  {\sf Einstein Manifolds} \cite{bes}, which     resulted in many of the  speculations discussed there  quickly becoming 
   folk-conjectures. 
   One such  speculation    \cite[p. 19]{bes}  concerned  whether the sign of the Einstein constant is somehow determined by the 
   diffeomorphism-type of the  Einstein  manifold; that is, if a smooth closed manifold $M$ admits two Einstein
   metrics $g_1$ and $g_2$, do the Einstein constants $\lambda_1$ and $\lambda_2$ of these two metrics necessarily  have the same sign? 
   However, contrary to Besse's apparent   expectations, the answer to this question was  later shown  to be {\em no}  by Fabrizio Catanese and
   the first author  \cite{cat}: 
   
  \begin{thm}[Catanese-LeBrun] \label{startled}
  For each $m\geq 2$, there is a closed simply-connected $4m$-manifold $M$ that carries both a $\lambda >0$ Einstein metric and a $\lambda < 0$ Einstein metric. 
  \end{thm}

A  frustrating limitation of the proof of this result, however,  was  that it only produced one example in each relevant dimension. This was because the entire construction depended on
  first producing a homeomorphic pair of simply-connected compact complex surfaces that were equipped with {\em K\"ahler}-Einstein metrics with Einstein  constants of 
  opposite sign. However, various breakthroughs in 
  complex surface theory eventually   made in possible for Rare\c{s} R\u{a}sdeaconu and Ioana \c{S}uvaina \cite{raresioana} to  improve Theorem \ref{startled} 
 in the following substantial manner:

  \begin{thm}[R\u{a}sdeaconu-\c{S}uvaina] 
  \label{improvement}
  For every $m \geq 2$, there are at least $\binom{m + 3}{3}$ distinct smooth closed simply-connected $4m$-manifolds
  that admit both $\lambda > 0$ and $\lambda < 0$ Einstein metrics. 
    \end{thm}

  Nonetheless, the technique used to prove Theorem \ref{improvement} was   based on roughly  the same strategy that  had previously been used to 
  prove Theorem \ref{startled}; namely, one begins by  identifying  minimal complex surfaces of general type that are homotopy-equivalent to 
  del Pezzo surfaces, and then tries to show that these candidates can actually be deformed into complex surfaces with  $c_1< 0$. 
  Not only did this make every small bit of progress in the subject seem as delicate as threading a needle, but it also {\em a priori} confined  
  the scope of these constructions to a small finite number of possible diffeo-types in each allowed  dimension. In this article, though, we will 
  finally  throw off these shackles by  adding a host of  newly-discovered non-K\"ahler building-blocks \cite{li-sun-cpt} to  our repertory. 
  This will allow us to prove the following:

  \begin{main} \label{banquet} 
   For each $m\geq 2$, there are {\em infinitely many} different smooth compact $4m$-manifolds that 
  admit both an Einstein metric with $\lambda > 0$ and an Einstein metric with $\lambda < 0$.  
  \end{main}

  Our proof of this result will actually be far simpler and less technical than the proofs of  Theorems \ref{startled} and \ref{improvement}. 
  Related arguments will also  allow us to dispense with a  peculiar feature of \cite{cat} and \cite{raresioana}, 
  which was that all of 
  the examples 
  constructed there  were  non-spin. 
  Here, we will   show  that 
  this pecularity was merely a methodological 
  artifact, by  displaying a sequence of  spin examples:
  
  \begin{main} \label{dessert} 
  \label{focus} 
  For each $m\geq 2$, there is a smooth closed {\em spin} $4m$-manifold  that 
  admits a  pair of Einstein whose Einstein constants have opposite signs. 
    \end{main}
  
  \section{Toric Einstein Manifolds} 
  
Our proof of Theorem \ref{banquet} will  make extensive use of   the following  surprising  result 
  \cite[Theorem 1.2] {li-sun-cpt}, discovered by  the second author and Song Sun:  
  
  \begin{thm}[Li-Sun] \label{abundant} The smooth compact $4$-manifold $M=  \mathsf{p} \CP_2 \# \mathsf{q} \overline{\CP}_2$
  admits a $T^2$-invariant Einstein metric iff its Euler characteristic $\chi (M)$ and signature $\tau (M)$ satisfy the inequality 
  \begin{equation}
  \label{balance} 
\chi (M) \geq 3 |\tau (M)|.
\end{equation}
  \end{thm}

  In the special case  where equality holds in \eqref{balance},  such ``toric'' Einstein metrics 
  were also discovered independently (and essentially simultaneously) by Tristan Ozuch 
  \cite{ozuch-toric}.  
  However, our proof of Theorem \ref{banquet} will 
  make essential use of  the  greater  flexibility that is  afforded by allowing strict inequality to hold  in \eqref{balance}.
Still, we find it plausible that an alternative proof of Theorem \ref{banquet} 
might   be  feasible 
that would only involve    examples that saturate  the 
inequality  \eqref{balance}. We certainly  hope 
 that some   readers of this article will decide  to take on  this technical challenge!

 The Einstein metrics of Theorem \ref{abundant} are said to be {\em toric}  because they admit effective isometric actions of the $2$-torus $T^2$. These Einstein
  metrics therefore  carry a pair of {\em Killing fields} 
 that generate the two relevant  circle actions on the relevant $4$-manifolds $M$.  Because each such  $4$-manifold $M$  is simply connected, 
 its  Euler characteristics $\chi (M)=2+b_2(M)$ is automatically  positive, and, by the Poincar\'e-Hopf  Theorem,  each  Killing field must therefore
 have  
 a zero\footnote{Indeed, by  \cite{kobfixed}, the $T^2$-action must actually have  $\chi (M)$ isolated fixed points.}  
somewhere on $M$. It therefore follows that Einstein constant $\lambda$ of any such toric Einstein $4$-manifold 
 must {\em a priori} be positive. Indeed, a so-called  ``Weitzenb\"ock'' argument due to Bochner \cite{bochner} shows  that 
 a compact manifold with non-positive  Ricci curvature can only carry a non-trivial Killing field if the field is parallel, 
 and so has constant length. This can then be used to show that, more generally,  any  compact  Einstein $4$-manifold that 
 carries a non-trivial Killing field  must either  have 
  $\lambda > 0$ or else be   finitely covered by a flat    $4$-torus. 
  
 For many choices of $\mathsf{p}$ and $\mathsf{q}$,   it  is also worth pointing out  that the proof of Theorem \ref{abundant} actually produces various  Einstein metrics 
on $M=\mathsf{p} \CP_2 \# \mathsf{q} \overline{\CP}_2$ that are invariant under different  $T^2$-actions on $M$. 
However, this interesting feature of \cite{li-sun-cpt} will not play an essential role in  the present article.

 The Einstein $4$-manifolds featured in  Theorem \ref{abundant} are all non-spin. However, 
 Li-Sun \cite{li-sun-cpt}  and Ozuch \cite{ozuch-toric} also constructed a  sequence of spin examples, 
and  other independent constructions of these   were then announced   shortly  afterwards by both Broder and Lafuente \cite{broder-lafuente} and by 
 Chu and Stern \cite{chu-stern-toric}: 
 
 \begin{thm} 
 \label{spin-dry} For each positive integer $\ell$, the  connected sum 
 $$\ell (S^2 \times S^2):= \underbrace{(S^2 \times S^2)\# \cdots \# (S^2 \times S^2)}_\ell$$
  admits a toric Einstein metric. 
 \end{thm} 
 
  This result will play a central role in our proof of Theorem \ref{dessert}.

  \section{K\"ahler-Einstein Partners} 
  
  \begin{prop} \label{partners} 
  Let $C_1,  \ldots , C_k$  be  pairwise  transverse  smooth cubic  hypersurfaces in $\CP_{k+2}$, 
  and  let $X_k=  C_1 \cap \cdots \cap C_k$ be the smooth compact complex surface in
  which they meet. If $k \geq 2$, then $X_k$ admits a $\lambda < 0$ K\"ahler-Einstein metric $g^-$ that is
  compatible with its complex structure. Moreover, the complex surface 
  $X_k$ is orientedly homeomorphic, but not diffeomorphic, 
  to the connected sum $\mathsf{p}_k \CP_2 \#\mathsf{q}_k \overline{\CP}_2$, where 
  \begin{eqnarray}
\mathsf{p}_k &=& 3^{k-1} (3k^2-7k+6)-1 , \label{mug}
\\
\mathsf{q}_k &=& 3^{k-1} (3k^2 + k+3) -1 .\label{wamp}
\end{eqnarray}
Consequently, 
$X_k$ is in particular $h$-cobordant to $\mathsf{p}_k \CP_2 \# \mathsf{q}_k \overline{\CP}_2$. 
  \end{prop}
  \begin{proof}
  By induction, the Lefschetz hyperplane  theorem implies that $X_k$ is simply connected. 
On the other hand,  since the anti-canonical bundle of  $\CP_{k+2}$ is isomorphic to $\mathcal{O}(k+3)$,  while the 
 normal bundle $\nu$ of  $X_k\subset \CP_{k+2}$ is isomorphic to the restriction of 
 $[\mathcal{O}(3)]^{\oplus k}$ to $X_k$,   adjunction implies that the canonical line  bundle $K$ of $X_k$ 
 is consequently  given by
 \begin{equation}
\label{connie}
K= \mathcal{O} (2k -3)|_{X_k}.
\end{equation}
In particular, if $\alpha \in H^2 (\CP_{k+2}, \ZZ)\cong \ZZ$ is the generator Poincar\'e dual to $\CP_{k+1}\subset \CP_{k+2}$, 
and if  $\iota : X_k\hookrightarrow \CP_{k+2}$ is the inclusion map, then \eqref{connie} implies that 
$c_1(X_k)= j^* (3-2k)\alpha$. 
 Since $X_k\subset \CP_{k+2}$ has degree $3^k$, it therefore follows    that 
 \begin{equation}
\label{seetoo}
c_1^2 (X_k) = 3^k (2k-3)^2.
\end{equation}
 Since the right-hand side here  is $\equiv 1 \bmod 2$,  this implies  that the
 intersection form of $X_k$ is odd, and  that the smooth   $4$-manifold 
 $X_k$ is consequently non-spin. Moreover, provided $k\geq 2$, equation \eqref{connie} also tells
 us that the canonical line bundle $K$ of $X_k$ is ample, so that 
 the Aubin-Yau theorem \cite{aubin,yau} therefore   implies that $X_k$ carries a unique   $\lambda = -1$ 
 K\"ahler-Einstein metric $g^-$ that is compatible with the constructed  complex structure on $X_k$.

Since   $p_1(T\CP_{k+2}) = (k+3) \alpha^2$, while the normal  bundle $\nu \to X_k$
has first Pontryagin class  $p_1(\nu) =\iota^*[ k (3\alpha)^2] = \iota^*[9 k \alpha^2]$,  adjunction therefore
yields
$$p_1(X_k) = p_1 (TX_k) = [(k+3) - 9k]\iota^* \alpha^2 =-  [8k-3]\iota^* \alpha^2,$$
and the Thom-Hirzebruch signature theorem therefore tells us that 
\begin{equation}
\label{siggi}
\tau (X_k) = \frac{1}{3} \int_{X_k}  p_1 (X_k) =-\frac{3^k}{3} (8k-3)= -3^{k-1} (8k-3).
\end{equation}
where $\tau(X_k)= b_+(X_k)-b_-(X_k)$ is by definition the signature of $X_k$. 
On the other hand, 
 \eqref{seetoo} and \eqref{siggi} imply that  the Euler characteristic of $X_k$ is given by 
\begin{eqnarray}
\chi (X_k) &=& \frac{c_1^2 (X_k) -3 \tau (X_k)}{2} \nonumber
\\ &=& \frac{3^k [(2k-3)^2 + (8k-3)]}{2} \nonumber
\\&=& 3^k(2k^2 - 2k +3).\label{oily} 
\end{eqnarray}
Since $X_k$ is simply connected, we therefore have 
\begin{eqnarray*}
1+ b_+(X_k) &=& \frac{\chi(X_k) + \tau (X_k)}{2}  \\
&=&  \frac{ 3^k(2k^2 - 2k +3) -3^{k-1} (8k-3) }{2}\\
&=& 3^{k-1} (3k^2 - 7k +6) 
\end{eqnarray*}
and 
\begin{eqnarray*}
1+ b_-(X_k) &=& \frac{\chi(X_k) - \tau (X_k)}{2}\\ &=& \frac{ 3^k(2k^2 - 2k +3) +3^{k-1} (8k-3) }{2}
\\&=&  3^{k-1} (3k^2 + k +3) ,
\end{eqnarray*}
so that, setting $\mathsf{p}_k:= b_+(X_k)$ and $\mathsf{q}_k:= b_-(X_k)$, these numbers are then given by
  \eqref{mug} and \eqref{wamp}.
Since Freedman's  results \cite{freedman} and 
the classification of indefinite quadratic forms \cite{hm}    show  
 that two smooth compact  simply-connected non-spin $4$-manifolds 
with indefinite intersection forms are orientedly homeomorphic iff their invariants $b_\pm$ agree, it therefore   follows that 
$X_k$ is homeomorphic to $\mathsf{p}_k\CP_2 \# \mathsf{q}_k \overline{\CP}_2$. Moreover, a classic  result of CTC Wall \cite{wall}
similarly  shows that $X_k$ and $\mathsf{p}_k\CP_2 \# \mathsf{q}_k \overline{\CP}_2$ are $h$-cobordant. On the other hand, 
for $k \geq 2$, $X_k$ and $\mathsf{p}_k\CP_2 \# \mathsf{q}_k \overline{\CP}_2$ are certainly not diffeomorphic, 
because Seiberg-Witten theory \cite{spccs,witten} implies that $X_k$ cannot admit metrics of positive scalar curvature, while
 $\mathsf{p}_k\CP_2 \# \mathsf{q}_k \overline{\CP}_2$ most certainly does admit such metrics, by \cite{gvln}.
  \end{proof}

  \begin{cor} \label{ranger} 
  The compact complex surface $X_k$ of Proposition \label{partners} satisfies
 $$\chi (X_k) \geq 3|\tau (X_k)|$$
 if and only if  $k\geq 5$. When this happens, the inequality is moreover strict. 
  \end{cor} 
  \begin{proof} By \eqref{siggi} and \eqref{oily}, 
  $$\chi(X_k) - 3|\tau (X_k)| = 3^k(2k^2 - 2k +3)-3\cdot 3^{k-1} (8k-3)= 2\cdot 3^k (k^2-5k+3).$$
  The claim thus follows from the fact that the asymptotically positive quadratic function 
  $x^2 - 5x +3$ has two real roots $(5\pm \sqrt{13})/2$, and that these two roots are respectively located in   the 
  open intervals $]0,1[$ and $]4,5[$. 
  \end{proof} 
  
  \section{Proofs of the Main Results}

  Theorem \ref{banquet} will now be an immediate corollary of the following: 
  
  \begin{thm} \label{fireworks} 
   For each integer $k\geq 5$,  once again consider the smooth compact simply-connected non-spin 
  $4$-manifold $M_k:= \mathsf{p}_k\CP_2 \# \mathsf{q}_k\overline{\CP}_2$, where 
   \begin{eqnarray*}
\mathsf{p}_k &=& 3^{k-1} (3k^2-7k+6)-1  \quad \text{and}
\\
\mathsf{q}_k &=& 3^{k-1} (3k^2 + k+3) -1 \end{eqnarray*}
are as    defined by  equations \eqref{mug} and \eqref{wamp}. Then, for 
  any integer $m \geq 2$, the smooth compact  non-spin $4m$-manifold
  $$M_k^{\times m} := \underbrace{M_k \times \cdots \times M_k}_m$$
  admits both a $\lambda > 0$  Einstein metric   and a $\lambda < 0$ Einstein metric. 
  \end{thm}
  \begin{proof} By Proposition \ref{partners}, the simply-connected complete-intersection complex surface $X_k$
  admits a  $\lambda = -1$ Einstein metric $g^-$, and is  $h$-cobordant to  $M_k$. 
  In particular, the $4m$-manifold 
  $$X_k^{\times m} := \underbrace{X_k \times \cdots \times X_k}_m$$
  therefore admits the the  $\lambda = -1$ Einstein metric $g_-:=g^-\oplus \cdots \oplus g^-$.
  But $X_k^{\times m}$ is also $h$-cobordant to $M_k^{\times m}$; indeed, if $W_k$ is  a compact 
 5-manifold-with-boundary such that 
 $$\partial W_k = \overline{X}_k \sqcup M_k$$
 and such that   $X_k \hookrightarrow W_k$ and $M_k \hookrightarrow W_k$  are both homotopy equivalences, then 
we can construct an $h$-cobordism from $X_k^{\times m}$ to $M_k^{\times m}$ by gluing the $m$ different $h$-cobordisms
$$\underbrace{X_k \times  \cdots \times X_k }_{m-j} \times W_k \times \underbrace{M_k \times \cdots \times M_k}_{j-1}, \qquad j= 1, \ldots , m$$
end-to-end. Thus, provided $m\geq 2$, Smale's $h$-cobordism theorem \cite{milnorh,smale} implies that  there is a 
diffeomorphism between the simply-connected $4m$-manifolds $X_k^{\times m}$  and $M_k^{\times m}$. Using this diffeomorphism to transplant 
$g^-$ from $X_k^{\times m}$ to $M_k^{\times m}$ then shows  that $M_k^{\times m}$  admits a $\lambda = -1$ Einstein metric.

On the other hand,  when 
  $k \geq 5$, Corollary \ref{ranger} and the existence of  a homotopy equivalence $M_k\sim X_k$ 
 together  tell us that $\chi (M_k ) > 3|\tau (M_k)|$, so 
   Theorem \ref{abundant}  guarantees that  the connected sum $M_k$ admits a 
  $\lambda = +1$ Einstein metric $g^+$. However,  this of course implies  that $M_k^{\times m}$ admits
  the $\lambda = +1$ Einstein metric $g_+:=g^+\oplus \cdots \oplus g^+$. This endows $M_k^{\times m}$ with 
  a pair of Einstein metrics  with Einstein constants of opposite signs, as promised.
   \end{proof}
  
  Since $\chi (M_k^{\times m}) = [\chi (M_k)]^m\to \infty$, the manifolds $M_k^{\times m}$ realize 
  infinitely many homotopy types, and Theorem \ref{banquet} is  therefore an immediate consequence. 
 Of course, for any $m$-tuple $(k_1, \ldots , k_m)$ of integers $k_j \geq 5$, where $m\geq 2$, 
 the same argument also shows that $X_{k_1} \times \cdots \times X_{k_m}$ and 
 $M_{k_1} \times \cdots \times M_{k_m}$  are diffeomorphic manifolds that admit Einstein metrics with
 Einstein constants of opposite signs. Moreover, many other  choices of the complex surfaces $X_k$ 
 would  certainly be possible. However, the particular choices that we have made here have certainly 
 sufficed to provide   a satisfyingly direct proof of Theorem \ref{banquet}, thereby realizing the  main goal of  this article. 
 
 \medskip

Our proof of Theorem \ref{dessert} is  less elementary and self-contained. However,  this result is in any case an  immediate corollary  of   the following: 
 
 \begin{thm} Let 
 $$Y = \underbrace{(S^2 \times S^2) \# \cdots \# (S^2 \times S^2)}_{56! (35)^2-1}.$$ 
 Then, for any $m\geq 2$,  the smooth compact  spin $4m$-manifold 
 $$Y^{\times m} :=\underbrace{Y \times  \cdots \times Y}_{m}$$
 admits a pair of  Einstein metrics with  Einstein constants of opposite signs.  \end{thm} 
 \begin{proof} 
 For $k=0$, $a=7$, and $b=4$, 
 we will  let $Z$ 
 be the  branched cover $\tilde{F}_{k, (a,b)}=\tilde{F}_{0, (7,4)}$ of 
 $\CP_1\times \CP_1$
 defined by Moisehzon, Robb,  and 
 Teicher \cite{morote}. Then  \cite[Theorem 0.4a]{morote} tells us that $Z$ is spin, simply connected,  and has signature $\tau (Z) = 0$. Since 
  \cite[Lemmas 7.1.1a and 7.1.3a,b]{morote}  also tell us that   $c_1^2 (Z) = n! (m-6)^2/4$, where $n= 2ab+kb^2 = 56$ and $m= 6ab-2a-2b+k(3b^2-b) = 146$,
  the topological Euler characteristic of $Z$ is given by $$\chi (Z) = \frac{c_1^2 (Z) - 3\tau (Z)}{2} = \frac{56! (140)^2}{8} = 2\cdot 56! (35)^2.$$
  Moreover,  the proof of \cite[Proposition 7.2.1]{morote} shows that $Z$ has ample canonical line bundle, 
  and it therefore follows that $Z$ admits a $\lambda= -1$ K\"ahler-Einstein metric $g^-$ by the Aubin-Yau theorem \cite{aubin,yau}. 
  Hence the product $Z^{\times m}=Z \times \cdots \times Z$ also admits the $\lambda=-1$ Einstein metric $g^-\oplus \cdots \oplus g^-$.

 But Wall's theorem \cite{wall} and the  classification of even unimodular quadratic forms \cite{hm}  also  tell us that
$Z$ is $h$-cobordant  to the connected sum $Y$.  Since this implies  that 
the simply-connected $4m$-manifold $Z^{\times m}$ and $Y^{\times m}$ 
are $h$-cobordant, the $h$-cobordism theorem \cite{milnorh,smale} therefore tells us that $Z^{\times m}$ is 
diffeomorphic to $Y^{\times m}$ whenever $m\geq 2$. This allows us to view $g^-\oplus \cdots \oplus g^-$
as a $\lambda =-1$ Einstein metric on $Y^{\times m}$. However, Theorem \ref{spin-dry}  also provides $Y$ with a
$\lambda = +1$ metric $g^+$, and $Y^{\times m}$ therefore also carries a $\lambda=+1$ Einstein metric $g^+\oplus \cdots \oplus g^+$, 
  as promised. 
 \end{proof} 
 
 Here, we note in passing that Moishezon, Robb, and Reicher actually  discovered  four different simply-connected spin surfaces with signature zero and ample canonical 
 line bundle in \cite{morote}, namely the branched cover $\tilde{F}_{k, (a,b)}$  of the $k^{\rm th}$ Hirzebruch surface $F_k$, for 
  $k = 0,1,2,3$,   $a= 7-2k$ and $b=4$. However, they do not seem to have noticed that these complex surfaces are all in fact  mutually homeomorphic! Thus, 
  one can actually   obtain $\binom{m+3}{3}$  different $\lambda =-1$ K\"ahler-Einstein metrics on $Y^{\times m}$ by essentially the same construction as described above. 
  It  would certainly be interesting to know whether or not these various Einstein metrics actually belong to different connected components of the Einstein moduli space. 
  It would also be natural and interesting  to compare and contrast this phenomenon with the situation  for $\lambda =+1$, where one  obtains many different 
  toric Einstein metrics on $Y^{\times m}$ that arise from  various different $T^2$ actions on   $Y$.
 
 \bigskip 
 
Finally, let us point out that, for any of the $4m$-manifolds used to prove Theorem \ref{banquet} and \ref{dessert},
the total volume of the $\lambda = +1$  Einstein metric is always {\em smaller} than that of the $\lambda = -1$ Einstein metric.
This happens because a theorem of Gursky \cite{G1,gl2} implies that  any positive-$\lambda$ Einstein metric on a 
compact oriented $4$-manifold $M$ with $b_+(M)\geq 2$ must satisfy the strict inequality $\lambda^2 V < 2\pi (2\chi +3\tau )(M)$,
where $V$ is the total volume, while, by contrast,  any K\"ahler-Einstein metric on a compact $4$-manifold must instead satisfy $\lambda^2 V = 2\pi (2\chi +3\tau )(M)$,
independent of the sign of $\lambda$. This should be compared  with the situation in \cite{cat,raresioana}, where, as a curious methodological  artifact,  the two Einstein metrics  
always 
had
identical total volumes.


\pagebreak 
%

 \vfill 
  
\noindent
{\bf Authors' addresses:}

\bigskip

\noindent
Claude LeBrun\\
Department of Mathematics\\
Stony Brook University\\
Stony Brook, NY 11794-3651

\bigskip

\noindent
Mingyang Li\\
Simons Center for Geometry and Physics\\
Stony Brook University\\
Stony Brook, NY 11794-3636

\end{document}